\documentclass{amsart}
\usepackage{palatino}
\usepackage{amsthm, amsmath}
\usepackage{amssymb} 
\usepackage{amscd}
\usepackage[hmargin=3cm, vmargin=3cm]{geometry}
\usepackage[breaklinks,bookmarksopen,bookmarksnumbered]{hyperref}
\usepackage[all]{xy}

\theoremstyle{plain}

\newtheorem{lem}{Lemma}
\newtheorem{prop}{Proposition}
\newtheorem{cor}{Corollary}

\theoremstyle{remark}
\newtheorem{rem}{Remark}
\newtheorem{ex}{Example}

\newcommand\pr{\noindent\textit{Proof} : }

\newcommand\rond{\kern 1pt{\scriptstyle\circ}\kern 1pt}

\def\lr#1{\langle {#1} \rangle}
\newcommand\End{\operatorname{End}}

\newcommand\Hom{\operatorname{Hom}}

\newcommand\rk{\operatorname{rk}}

\newcommand\Z{\mathbb{Z}}
\newcommand\Q{\mathbb{Q}}
\newcommand\R{\mathbb{R}}
\newcommand\C{\mathbb{C}}
\renewcommand\P{\mathbb{P}}

\newcommand\F{\mathbb{F}}

\renewcommand\O{\mathcal{O}}

\newcommand\NS{\mathrm{NS}}

\def\qfl#1{\buildrel {#1}\over {\longrightarrow}}
\newcommand\mono{\lhook\joinrel\mathrel{\longrightarrow}}
\newcommand\iso{\vbox{\hbox to .8cm{\hfill{$\scriptstyle\sim$}\hfill}
\nointerlineskip\hbox to .8cm{{\hfill$\longrightarrow $\hfill}} }}

\begin{document}
\title[Some surfaces with maximal Picard number]{Some surfaces with maximal Picard number}
\author[Arnaud Beauville]{Arnaud Beauville}
\address{Laboratoire J.-A. Dieudonn\'e\\
UMR 7351 du CNRS\\
Universit\'e de Nice\\
Parc Valrose\\
F-06108 Nice cedex 2, France}
\email{arnaud.beauville@unice.fr}
 
 
\begin{abstract}
For a smooth complex projective variety, the rank $\rho $ of the N\'eron-Severi group is bounded by the Hodge number $h^{1,1}$. Varieties with $\rho =h^{1,1}$ have interesting properties, but are rather sparse, particularly in dimension 2. We discuss in this note a number of examples, in particular those constructed from curves with special Jacobians. 
\end{abstract}

\maketitle 
\section{Introduction}

The \emph{Picard number} of  a smooth projective variety $X$ is the rank $\rho $ of the N\'eron-Severi group -- that is, the group of  classes of divisors in $H^2(X,\Z)$. It is bounded by the Hodge number $h^{1,1}:=\allowbreak \dim H^1(X,\Omega ^1_X)$. We are interested here in varieties with maximal Picard number $\rho =h^{1,1}$. As we will see in \S 2, there are many examples of such varieties in dimension $\geq 3$,  so we will focus on the case of surfaces. 

Apart from the well understood case of K3 and abelian surfaces, the quantity of known examples is remarkably small. In \cite{P} Persson showed that some families of
 double coverings of rational surfaces contain surfaces with maximal Picard number (see section \ref{P} below);   some scattered examples have appeared since then. We will review them in this note  and examine in particular when the product of two curves has maximal Picard number -- this provides some examples, unfortunately also quite sparse.
 
 \medskip	
\section{Generalities}
Let $ X $ be a smooth projective variety over $ \C $. The N\'eron-Severi group $ \NS(X)  $ is the subgroup of algebraic classes in $ H^2(X,\Z)  $; its rank $ \rho $ is the \emph{Picard number} of $ X $. The natural map $ \NS(X) \otimes\C\rightarrow H^2(X,\C)  $ is injective and its image is contained in $ H^{1,1} $, hence  $ \rho\leq h^{1,1} $.

\begin{prop}\label{triv}
The following conditions are equivalent :

{\rm (i) }  $ \rho= h^{1,1} $;

{\rm (ii) } The  map $ \NS(X) \otimes\C\rightarrow H^{1,1}   $ is bijective;

{\rm (iii) } The subspace $ H^{1,1}$ of $ H^2(X,\C) $ is defined over $ \Q $. 

{\rm (iv) }  The subspace $ H^{2,0}\oplus H^{0,2} $ of $ H^2(X,\C) $ is defined over $ \Q $. 
\end{prop}
\pr The  equivalence of {\rm (iii)} and {\rm (iv)} follows from the fact that  $ H^{2,0}\oplus H^{0,2} $ is the orthogonal of $ H^{1,1} $ for the scalar product on $ H^2(X,\C) $ associated to an ample class. The rest is clear.\qed

\medskip	
When $X$ satisfies these equivalent properties we will say for short that $X$ is
\emph{$ \rho $-maximal} (one finds the terms singular, exceptional or extremal in the literature).

\medskip	
\noindent\emph{Remarks}$ .- $
1) A variety with   $ H^{2,0}=0 $ is  $ \rho $-maximal. We will implicitly exclude this trivial case in the discussion below.

\smallskip	
2) Let $X$, $Y$ be two $\rho $-maximal varieties, with $H^1(Y,\C)=0$. Then $X\times Y$ is $\rho $-maximal. For instance $X\times \P^n$ is $\rho $-maximal, and $Y\times C$ is $\rho $-maximal for any curve $C$.  

\smallskip	
3) Let $ Y $ be a submanifold of $ X $; if $ X $ is $ \rho $-maximal and the restriction map $ H^2(X,\C)\rightarrow  H^2(Y,\C) $ is bijective,  $ Y $ is $ \rho $-maximal. By the Lefschetz theorem, the latter condition is realized if $ Y $ is a complete intersection of smooth ample divisors in $ X $, of dimension $ \geq 3 $. Together with Remark 2, this gives many examples of  $ \rho $-maximal varieties of dimension $ \geq 3 $; thus we will focus on finding $ \rho $-maximal \emph{surfaces}.

\setcounter{rem}{3}
\begin{prop}\label{map}
Let $\pi: X\dasharrow Y $ be a rational map of smooth projective varieties. 

$ a) $ If  $\pi ^*:H^{2,0}(Y)\rightarrow H^{2,0}(X)$ is injective (in particular if $\pi $ is dominant), and
$ X $ is $ \rho $-maximal, so is $ Y $.

$ b) $ If $\pi ^*:H^{2,0}(Y)\rightarrow H^{2,0}(X)$ is surjective and $ Y $ is $ \rho $-maximal, so is $ X $.
\end{prop}
Note that since $\pi $ is defined on an open subset $U\subset X$ with $\mathrm{codim}(X\smallsetminus U)\geq 2$, the pull back map $\pi ^*: H^2(Y,\C)\rightarrow H^2(U,\C)\cong H^2(X,\C)$ is well defined.

\smallskip	
\pr  Hironaka's theorem provides a diagram
\[ \xymatrix{& {\hat X} \ar[dl]_b \ar[dr]^{\hat \pi} & \\
X\ar@{-->}[rr]^\pi & & Y} \]
where $ \hat{\pi} $ is a morphism, and $ b $ is a composition of blowing-ups with smooth centers. Then  $b^*: \allowbreak H^{2,0}(X)\rightarrow H^{2,0}(\hat{X})$ is bijective,  and
$\hat{ X} $ is $ \rho $-maximal if and only if $ X $ is $ \rho $-maximal; so replacing $ \pi $ by $ \hat{\pi} $ we may assume
 that $ \pi $ is a morphism.

$a)$ Let $V:=(\pi ^*)^{-1}(\mathrm{NS}(X)\otimes \Q)$. We have 
  \[V\otimes_{\Q} \C= (\pi ^*)^{-1}(\mathrm{NS}(X)\otimes \C)=(\pi ^*)^{-1}(H^{1,1}(X))=H^{1,1}(Y)\]
  (the last equality holds because
  $\pi ^*$ is injective on $H^{2,0}(Y)$ and $H^{0,2}(Y)$), hence $Y$   is $\rho $-maximal.

 $b)$  Let $W$ be the $\Q$-vector subspace of $H^2(Y,\Q)$ such that $W\otimes _{\Q}\C= H^{2,0}(Y)\oplus H^{0,2}(Y)$. Then $\pi ^*W$ is a $\Q$-vector subspace of $H^2(X,\Q)$, and 
 \[(\pi ^*W)\otimes \C = \pi ^*(W\otimes \C)=\pi ^*H^{2,0}(Y)\oplus \pi ^*H^{0,2}(Y)=H^{2,0}(X)\oplus H^{0,2}(X)\]
 so $X$ is $\rho $-maximal.
\qed

\medskip
\section{Abelian varieties}
There is a nice characterization of $\rho $-maximal abelian varieties   (\cite{K}, \cite{L}) :
\begin{prop}\label{ab}
Let $ A $ be an abelian variety of dimension $ g $. We have 
$  \rk_{\Z} \End(A) \leq 2g^2$.
The following conditions are equivalent :

{\rm(i) } $ A $ is $ \rho $-maximal;

{\rm(ii) }  $ \rk_{\Z} \End(A)=2g^2 $;

{\rm(iii) } $ A $ is isogenous to $ E^g $, where $ E $ is an elliptic curve with complex multiplication.

{\rm (iv) } $A$ is isomorphic to a product of mutually isogenous elliptic curves with complex multiplication.

\end{prop}

(The equivalence of (i), (ii) and (iii) follows easily from Lemma \ref{hom} below; the only delicate point is (iii)$\ \Rightarrow\ $(iv), which we will not use.)

\medskip	
Coming back to the surface case,
suppose that our abelian variety $ A $ contains a surface $ S $ such that the restriction map $ H^{2,0}(A)\rightarrow H^{2,0}(S)   $ is surjective. Then $  S$ is $ \rho $-maximal    if $ A $ is $ \rho $-maximal (Proposition \ref{map}.$b$). Unfortunately this situation seems to be rather rare. We will discuss below (Proposition \ref{c}) the case of $\mathrm{Sym}^2C$ for a curve $C$. Another interesting example is the \emph{Fano surface} $ F_X $ parametrizing the lines contained in a smooth cubic threefold $ X $, embedded in the intermediate Jacobian $ JX $ \cite{CG}. There are some cases in which $ JX $ is known to be  $ \rho $-maximal :

\begin{prop}
$ a) $ For $ \lambda\in\C $, $ \lambda^3\neq 1 $, let $ X_{\lambda} $ (resp. $ E_{\lambda} $) be the cubic in $ \P^4 $ (resp. $\P^2$) defined by
\[ X_{\lambda}: X^3+Y^3+Z^3-3\lambda XYZ +T^3+U^3=0\qquad 
 E_{\lambda}: X^3+Y^3+Z^3-3\lambda XYZ =0 \ .\]
 If $E_{\lambda}$ is isogenous to $E_0$, $ JX_{\lambda} $ and $ F_{X_{\lambda}} $   are $ \rho $-maximal. The set of  $\lambda\in\C $ for which this happens is  countably infinite. 
 
 \smallskip	
 $ b) $
Let $ X \subset \P^4$ be the Klein cubic threefold $\displaystyle \sum_{i\in \Z/5}X_i^2X_{i+1}=0 $. Then $ JX $ and $ F_X $ are $ \rho $-maximal.
\end{prop}
\pr Part 1) is due to Roulleau \cite{R} , who proves that $ JX_{\lambda} $ (for any $ \lambda $) is isogenous to  $ E_0^3\times E_{\lambda}^2 $. Since the family $(E_\lambda )_{\lambda \in \C}$ is not constant, there is a countably infinite set of  $\lambda\in\C $ for which $E_{\lambda}$ is isogenous to $E_0$, hence $ JX_{\lambda} $ and therefore  $ F_{X_{\lambda}} $   are $ \rho $-maximal.

\smallskip

Part 2) follows from a result of Adler \cite{A}, who proves  that $ JX $ is isogenous (actually isomorphic)  to $ E^5 $, where $ E $ is the elliptic curve whose endomorphism ring is the ring of integers of $ \Q(\sqrt{-11} )  $ (see also \cite{Rk} for a precise description of the group $ \NS(X) $).\qed

\medskip
\section{Products of curves}\label{prod}

\begin{prop}\label{cc'}
Let $C,C'$ be two smooth projective curves, of genus $g$ and $g'$ respectively. The following conditions are equivalent :

{\rm (i)} The surface $C\times C'$ is $\rho $-maximal;

{\rm (ii)} There exists an elliptic curve $E$ with complex multiplication such that $JC$ is isogenous to $E^g$ and $JC'$ to $E^{g'}$.

\end{prop}
\pr Let  $ p,p' $ be the projections from $ C\times C' $ to $ C $ and $C'$. We have 
\[ H^{1,1}(C\times C')= p^*H^2(C,\C)\oplus p'^*H^2(C',\C)\, \oplus\,\left( p^*H^{1,0}(C) \otimes p'^*H^{0,1}(C')\right) \, \oplus\, \left(p^*H^{0,1}(C) \otimes p'^*H^{1,0}(C') \right) \]
hence $ h^{1,1}(C\times C')=2gg'+2 $. On the other hand we have \[ \NS(C\times C') = p^*\NS(C)\,\oplus\,  p'^*\NS(C')\,\oplus\, \Hom (JC, JC')    \] (\cite{BL}, thm.\ 11.5.1), hence $ C\times C' $ is $ \rho $-maximal if and only if $ \rk \Hom (JC, JC')  = 2gg'  $. Thus the Proposition follows from the following (well-known) lemma :

\begin{lem}\label{hom}
Let $A$ and $B$ be two abelian varieties, of dimension $a$ and $b$ respectively. The $\Z$-module $\Hom(A,B)$ has rank $\leq 2ab$;  equality holds if and only if there exists an elliptic curve $E$ with complex multiplication such that $A$ is isogenous to $E^a$ and $B$ to $E^b$.
\end{lem}
\pr There exist simple abelian varieties $A_1,\ldots ,A_s$, with distinct isogeny classes, and nonnegative integers $p_1,\ldots ,p_s,q_1,\ldots ,q_s$ such that $A$ is isogenous to $A_1^{p_1}\times \ldots \times A_s^{p_s}$ and $B$ to $A_1^{q_1}\times \ldots \times A_s^{q_s}$. Then
\[\Hom(A,B)\otimes_{\Z}\Q\,\cong\, M_{p_1,q_1}(K_1)\times \ldots \times  M_{p_s,q_s}(K_s)\]
where $K_i$ is the (possibly skew) field $\End(A_i)\otimes _{\Z}\Q$. Put $a_i:=\dim A_i$. Since $K_i$ acts on $H^1(A_i,\Q)$ we have $\dim_{\Q}K_i\leq \allowbreak b_1(A_i)=2a_i$, hence
\[\rk \Hom(A,B)\leq \sum_i 2p_iq_ia_i \leq 2 (\sum p_ia_i)(\sum q_ia_i)=2ab\ .\]
The last inequality is strict unless $s=a_1=1$,  in which  case
 the first one is strict unless $\dim_{\Q}K_1=2$. The lemma, and therefore the Proposition, follow.\qed

\medskip	
The most interesting case occurs when $C=C'$. Then :
\begin{prop}\label{c}
Let $ C $ be a smooth projective curve. The following conditions are equivalent :

{\rm (i) } The Jacobian $ JC $ is $ \rho $-maximal;

{\rm (ii) } The surface $ C\times C $ is $ \rho $-maximal;

{\rm (iii) } The symmetric square $\mathrm{Sym}^2C$ is $ \rho $-maximal.
\end{prop}

\pr  The equivalence of (i) and  (ii) follows from Proposition \ref{cc'}. The Abel-Jacobi map $\mathrm{Sym}^2C\rightarrow JC$ induces an isomorphism $ H^{2,0}(JC)\cong \wedge^2 H^{0}(C,K_C)\iso H^{2,0}(\mathrm{Sym}^2C)$, thus (i) and (iii) are equivalent by Proposition \ref{map}.\qed

\medskip	
When the equivalent conditions of Proposition \ref{c} hold, we will say that  $ C $ has \emph{maximal correspondences} (the group $ \End(JC)  $ is often called the group of divisorial correspondences of $ C $). 

By Proposition \ref{ab} the Jacobian $JC$ is then isomorphic to a product of elliptic curves
(\emph{completely decomposable in the strict sense} in the terminology of \cite{ES}). Though
we know very few examples of such curves, we will give below some examples with $g=4$ or $10$; this answers positively Question 3 of \emph{loc. cit.}

  For $ g=2 $ or $ 3 $, there is a countably infinite set of  curves with maximal correspondences (\cite{HN},\cite{H}). The point is that  any indecomposable principally polarized abelian variety of dimension 2 or 3 is a Jacobian; thus it suffices to construct an  indecomposable principal polarization on $ E^g $, where $ E $ is an elliptic curve with complex multiplication, and this is easily translated into a problem about hermitian forms of rank $ g $ on certain rings of quadratic integers. 
  
  This approach works only for $g=2$ or $3$; moreover it does not give an explicit description of the curves. Another method  is by using automorphism groups, with the help of the following easy lemma :
  
  \begin{lem}\label{irred}
Let $ G $ be a finite group of automorphisms of $ C $, and let $ H^0(C,K_C)=\oplus_{i\in I}V_i  $ be a decomposition of the $ G $-module $ H^0(C,K_C) $ into irreducible representations. Assume that there exists an elliptic curve $ E $ and for each $ i\in I $, a nontrivial map $ \pi_i:C\rightarrow E $ such that $ \pi_i^*H^0(E,K_E)\subset V_i  $. Then $ JC $ is isogenous to $ E^g $.
\end{lem} 
In particular if $H^0(C,K_C)$ is an irreducible $G$-module and $C$ admits a map onto an elliptic curve $E$, then $JC$ is isogenous to $ E^g $.

\medskip	
\pr Let $ \eta  $ be a generator of   $ H^0(E,K_E) $. Let $ i\in I $; the forms 
$ g^*\pi_i^*\eta  $ for $ g\in G$ generate $V_i$, hence there 
 exists a subset $ A_i $ of $ G $ such that the forms $ g^*\pi_i^*\eta $ for $ g\in A_i $ form a basis of $ V_i $.  

Put $ \Pi_i=(g\rond \pi_i)_{g\in A_i} :C\rightarrow E^{A_i} $, and  $ \Pi=(\Pi_i)_{i\in I}: C\rightarrow E^g  $. By construction $ \Pi^*:H^0(E^g,\Omega^1_{E^g})\rightarrow \allowbreak H^0(C,K_C)   $ is an isomorphism. Therefore the map $ JC\rightarrow E^g $ deduced from $ \Pi $ is an isogeny.\qed 

 \medskip	
 In the examples which follow, and in the rest of the paper, we put  $\omega :=e^{\frac{2\pi i}{3} }$.
\begin{ex} We consider the family $(C_t)$ of genus 2 curves given by $y^2=x^6+tx^3+1$, for $t\in \C\smallsetminus \{\pm 2\}$. It admits the automorphisms
$\tau :(x,y)\mapsto (\dfrac{1}{x}, \dfrac{y}{x^3 }) $ and $\psi :(x,y)\mapsto (\omega x,y)$.
  The forms $\dfrac{dx}{y}$ and $\dfrac{xdx}{y}  $ are eigenvectors for $\psi $ and  are exchanged (up to sign) by $\tau $; it follows that the action of the group generated by $\psi $ and $\tau $ on $H^0(C_t,K_{C_t})$ is irreducible. 
  
  Let $E_t$ be the elliptic curve defined by $v^2=(u+2)(u^3-3u+t)$; the curve $C_t$ maps onto $E_t$ by 
  $(x,y)\mapsto\left(x+\dfrac{1}{x}, \dfrac{y(x+1)}{x^2}  \right)$. By Lemma \ref{irred}  $JC_t$ is isogenous to $E_t^2$. 
 Since the $j$-invariant of $E_t$ is a non-constant function of $t$,  there is a countably infinite  set of $t\in\C$  for which $E_t$ has complex multiplication, hence $C_t$ has maximal correspondences.
\end{ex}

\begin{ex}
Let $ C $ be the genus 2 curve $ y^2=x(x^4-1)  $; its automorphism group is a central extension of $\mathfrak{S}_4$ by the hyperelliptic involution $\sigma $ (\cite{BL}, 11.7); its action  on $H^0(C,K_C)$ is irreducible. 

Let $E$ be the elliptic curve $E:\ v^2=u(u+1)(u-2\alpha )$, with $\alpha =1-\sqrt{2}$. 
The curve $C$ maps to $E$ by $(x,y)\mapsto \left(\dfrac{x^2+1}{x-1}, \dfrac{y(x-\alpha )}{(x-1 )^2}  \right)$. The 
 $j$-invariant of $E$ is $8000$, so
 $E$ is the elliptic curve $\C/\Z[\sqrt{-2}]$
 (\cite{S}, Prop. 2.3.1). 
\end{ex}

  \begin{ex}[The $\mathfrak{S}_4$-invariant quartic curves]\label{fermat}
Consider the standard representation of  $\mathfrak{S}_4$ on $\C^3$. It is convenient to view $\mathfrak{S}_4$ as the semi-direct product $(\Z/2)^2\rtimes\mathfrak{S}_3$, with $\mathfrak{S}_3$ (resp. $(\Z/2)^2$) acting on $\C^3$ by permutation (resp.\ change of sign) of the basis vectors. The quartic forms invariant under this representation form the pencil
\[(C_t)_{t\in\P^1} \ :\ x^4+y^4+z^4+t(x^2y^2+y^2z^2+z^2x^2)=0\ .\]
According to \cite{DK}, this pencil was known to Ciani. It contains the Fermat quartic $(t=0)$ and the Klein quartic $(t=\frac{3}{2}(1\pm i\sqrt{7}) )$.

Let us take $t\notin\{2,-1,-2,\infty\}$; then
$C_t$ is smooth. The action of $\mathfrak{S}_4$ on $H^0(C_t,K)$, given by the standard representation, is irreducible. Moreover the involution $x\mapsto -x$ has 4 fixed points, hence the quotient curve $E_t$ has genus 1. It is given by the degree 4 equation
\[u^2+tu (y^2+z^2) +y^4+z^4+ty^2z^2=0 \]
in the weighted projective space $\P(2,1,1)$. Thus $E_t$ is a double covering of $\P^1$ branched along the zeroes  of the polynomial $(t+2)(y^4+z^4)+2ty^2z^2$. The cross-ratio of these zeroes is $-(t+1)$, so $E_t$ is the elliptic curve $y^2= x(x-1)(x+t+1)$. By Lemma \ref{irred} $JC_t$ \emph{is isogenous to} $E_t^3$. For a countably infinite set of  $t$ the curve $E_t$ has complex multiplication, thus $C_t$ has maximal correspondences. 
For $t=0$ we recover the well known fact that the Jacobian of the Fermat quartic curve is isogenous to  
$(\C/\Z[ i])^3$.
\end{ex}

\begin{ex}
Consider the genus 3 hyperelliptic curve $H:\ y^2=x(x^6+1)$. The space $H^0(H,K_H)$ is spanned by $\dfrac{dx}{y},\, x\dfrac{dx}{y},\, x^2\dfrac{dx}{y}$. This is a basis of eigenvectors for the automorphism $\tau :(x,y)\mapsto (\omega x,\omega^2 y)$. On the other hand the involution $\sigma :(x,y)\mapsto (\dfrac{1}{x},\dfrac{-y}{x^4}  )$ exchanges $\dfrac{dx}{y}$ and $x^2\dfrac{dx}{y}$, hence the summands of the decomposition
\[H^0(H,K_H) = \lr{\frac{dx}{y}, x^2\frac{dx}{y}}\oplus \lr{x\frac{dx}{y}}\]
are irreducible under the group $\mathfrak{S}_3$ generated by $\sigma $ and $\tau $.

Let $E_i$ be the elliptic curve $v^2=u^3+u$, with endomorphism ring $\Z[i]$. Consider the maps $f$ and $g$ from $H$ to $E_i$ given by
\[f(x,y)=(x^2,xy) \qquad g(x,y)=(\lambda^2(x+\frac{1}{x})\ ,\ \frac{\lambda ^3y}{x^2}  )\quad\mbox{with}\ \lambda ^{-4}=-3\ . \]
We have $f^*\dfrac{du}{v}=\dfrac{2xdx}{y}  $ and $g^*\dfrac{du}{v}=\lambda^{-1} (x^2-1)\dfrac{dx}{y}$. Thus we can apply Lemma \ref{irred}, and we find that $JH$ \emph{is isogenous to $E_i^3$}. 

Thus $JH$ is isogenous to the Jacobian of
 the Fermat quartic $F_4$  (example \ref{fermat}). In particular we see that the surface $H\times F_4$ is $\rho $-maximal.

\end{ex}

  We now arrive to our main example in higher genus. Recall that we put $\omega = e^{\frac{2\pi i}{3}}$.
  
\begin{prop}\label{C6}
The Fermat sextic curve $ C_6: X^6+Y^6+Z^6=0 $ has maximal correspondences. Its Jacobian $JC_6$ is isogenous to $E_\omega ^{10}$, where $E_\omega $ is the elliptic curve $\C/\Z[\omega ]$.
\end{prop}

 The first part can be deduced from the general recipe given by Shioda to compute the Picard number of $ C_d\times C_d $ for any $ d $ \cite{S2}. Let us give an elementary proof. Let $G:= T\rtimes \mathfrak{S}_3$, where $\mathfrak{S}_3$ acts on $\C^3$ by permutation of the coordinates and 
 $T$ is the group of diagonal matrices $t$ with $t^6=1$. 
\[\Omega =\frac{XdY-YdX}{Z^5}=\frac{Y dZ-ZdY}{X^5}  = \frac{ZdX-XdZ}{Y^5}\ \in H^0(C,K_C(-3)) \ .\leqno{\mbox{\quad Let}}\]
A basis of eigenvectors for the action of $T$ on $H^0(C_6,K)$ is given by the forms $X^aY^bZ^c\,\Omega $, with $a+b+c=3$; using the action of $\mathfrak{S}_3$ we get a decomposition into irreducible components: 
\[H^0(C_6,K)=   V_{3,0,0}\oplus V_{2,1,0}\oplus V_{1,1,1}\]
where $V_{\alpha ,\beta ,\gamma }$ is spanned by the forms $X^aY^bZ^c\,\Omega$ with $\{a,b,c\}=\{\alpha ,\beta ,\gamma \}$.

Let us use affine coordinates $x=\frac{X}{Z} $, $y=\frac{Y}{Z} $ on $C_6$. We consider the following maps from $C_6$ onto $E_\omega :\ v^2=u^3-1$ :
\[f(x,y)= (-x^2,y^3) \qquad g(x,y)=\bigl(2^{-\frac{2}{3} }x^{-2}y^4\,,\, \frac{1}{2}(x^3-x^{-3})\bigr) \ ;\]
and, using for $E_{\omega }$ the equation  $\xi ^3+\eta ^3+1=0$, $h(x,y)=(x^2,y^2)$.

We have 
\[f^*\frac{du}{v}=-\frac{2xdx}{y^3}=-2 XY^2\,\Omega \in V_{2,1,0}  \quad,\quad g^*\frac{du}{v}=-2^{\frac{4}{3} }Y^3\,\Omega \in V_{3,0,0} \quad,\quad h^*\frac{d\xi }{\eta ^2}= 2XYZ\,\Omega \in V_{1,1,1}\] so the Proposition follows from Lemma \ref{irred}.\qed

\bigskip

By Proposition \ref{map} every quotient of $ C_6 $ has again maximal correspondences. There are four such quotient  which have genus 4 : 

$ \bullet $ The quotient by an involution  $\alpha \in T$, which we may take to be $ \alpha : (X,Y,Z)\mapsto (X,Y,-Z)  $. The canonical model of $ C_6/\alpha  $
is the image of $ C_6 $ by the map $ (X,Y,Z)\mapsto (X^2,XY,Y^2, Z^2)   $; 
its equations in $ \P^3 $ are $ xz-y^2= x^3+z^3+t^3=0 $. Projecting onto the conic $xz-y^2=0$ realizes $ C_6/\alpha  $ as the cyclic triple covering $v^3=u^6+1$ of $\P^1$.

$\bullet$ The quotient by an involution $\beta \in\mathfrak{S}_3$, say $\beta  :(X,Y,Z)\mapsto (Y,X,Z)$. The canonical model of $ C_6/\beta   $
is the image of $ C_6 $ by the map $ (X,Y,Z)\mapsto ((X+Y)^2, Z(X+Y), Z^2,XY)   $; its equations are $xz-y^2=x(x-3t)^2+z^3-2t^3=0$. 

Since the quadric containing their canonical model is singular, the two genus 4 curves $C_6/\alpha  $ and $C_6/\beta  $ have a unique $g^1_3$. The associated triple covering $C_6/\alpha  \rightarrow \P^1$ is  cyclic, while the corresponding covering $C_6/\beta  \rightarrow \P^1$
is not. Therefore the two curves are not isomorphic.

$ \bullet $ The quotient by an element  of order 3 of $T$ acting freely, say $ \gamma  : (X,Y,Z)\mapsto (X,\omega Y,\omega ^2Z)   $. The canonical model of $ C_6/\gamma  $  is the image of $ C_6 $ by the map $ (X,Y,Z)\mapsto (X^3,Y^3,Z^3,XYZ)$; its equations are 
$x^2+y^2+z^2=t^3-xyz=0$. Projecting onto the conic $x^2+y^2+z^2=0$ realizes $ C_6/\gamma   $ as the cyclic triple covering $v^3=u(u^4-1)$ of $\P^1$; thus
$C_6/\gamma $ is not isomorphic to $C_6/\alpha  $ or  $C_6/\beta  $.

$ \bullet $ The quotient by an element  of order 3 of $\mathfrak{S}_3$ acting freely, say $ \delta : (X,Y,Z)\mapsto (Y,Z,X)   $. The canonical model of $ C_6/\delta  $  is the image of $ C_6 $ by the map \[ (X,Y,Z)\mapsto (X^3+Y^3+Z^3,XYZ,X^2Y+Y^2Z+Z^2X, XY^2+YZ^2+ZX^2)   \ .\] It is  contained in the smooth quadric $ (x+y)^2+5y^2-2zt=0 $, so  $ C_6/\delta  $ is not isomorphic to any of the 3 previous curves. 

\medskip	
Thus we have found four non-isomorphic curves of genus 4 with Jacobian isogenous to $E_\omega ^4$. The product of any two of these curves is a $\rho $-maximal surface.

\begin{cor}\label{S6}
The Fermat sextic surface $ S_6: X^6+Y^6+Z^6+T^6=0$ is $ \rho $-maximal.
\end{cor}
\pr This follows from Propositions \ref{C6}, \ref{map} and Shioda's trick : there exists a rational dominant map $ \pi: C_6\times C_6 \dasharrow S_6 $, given by  $ \pi  \bigl((X,Y,Z),(X',Y',Z')\bigr)= (XZ',YZ',iX'Z,iY'Z) $.\qed

\medskip	
\begin{rem}
Since the Fermat plane quartic has maximal correspondences (Example 2), the same argument 
gives the classical fact that the Fermat quartic surface is  $ \rho $-maximal.
\end{rem}
\medskip	
Again every quotient of the Fermat sextic is $ \rho $-maximal. For instance, the quotient of $ S_6 $ by the automorphism $ (X,Y,Z,T)\mapsto  (X,Y,Z,\omega T) $ is the double covering of $ \P^2 $ branched along $ C_6 $: it is a $ \rho $-maximal K3 surface. The quotient of $S_6$ by the involution $(X,Y,Z,T)\mapsto (X,Y,-Z,-T)$ is given in $\P^5$ by the equations
\[y^2-xz = v^2-uw= x^3+z^3+u^3+w^3=0\ ;\]it is a complete intersection of degrees $(2,2,3)$, with 12 ordinary nodes. 
Other quotients have $p_g $ equal to 2,3, 4 or 6.

\medskip
\section{Quotients of self-products of curves}\label{cub}
The method of the previous section may sometimes allow to prove that certain quotients of a product $C\times C$ have maximal Picard number. 
Since we have very few  examples we will refrain from giving a general statement and contend ourselves with one significant example.

Let $C$ be the  curve in $\P^4$ defined by
\[ u^2=xy \quad,\quad v^2=x^2-y^2\quad,\quad w^2=x^2+y^2\ . \]
It is isomorphic to the modular curve $X(8)$ \cite{FS}. Let $\Gamma\subset \mathrm{PGL(5,\C)}$ be the subgroup of diagonal elements changing an even number of signs of $u,v,w$; $\Gamma $ is isomorphic to $(\Z/2)^2$ and acts freely on $C$.

 \begin{prop}
$a)$ $JC$  is isogenous to $E_i^3\times E^2_{\sqrt{-2}}$, where $E_\alpha =\C/\Z[\alpha ]$ for $\alpha =i$ or  ${\sqrt{-2}}$.

$b)$ The surface $(C\times C)/\Gamma $ is $\rho $-maximal.
\end{prop}

\pr $a)$ The form  $\ \Omega  :=(xdy-ydx)/uvw\ $ generates $H^0(C,K_C(-1))$, and is $\Gamma $-invariant; thus multiplication by $\Omega  $ induces a $\Gamma $-equivariant isomorphism $H^0(\P^4,\O_{\P^4}(1))\iso H^0(C,K_C)$. Let $V$ and $L$ be the subspaces of $H^0(C,K_C)$ corresponding to $\lr{u,v,w}$ and $\lr{x,y}$. The projection $(u,v,w,x,y)\mapsto \allowbreak  (u,v,w)$ 
maps $C$ onto the quartic curve $F:\ 4u^4+v^4-w^4=0$; the induced map $f:C\rightarrow F$
 identifies $F$ with the quotient of $C$ by the involution $(u,v,w,x,y)\mapsto(u,v,w,-x,-y)$, and we have $f^*H^0(F,K_F)=V$.
 
 The quotient curve $H:=C/\Gamma $ is the genus 2 curve $z^2=t(t^4-1)$ \cite{Bt}. The pull-back of $H^0(H,K_H)$ is the subspace invariant under $\Gamma $, that is $L$. Thus $JC$ is isogenous to $JF\times JH$. From examples  1 and 2 of \S \ref{prod} we conclude that $JC$ is isogenous to $E_i^3\times E^2_{\sqrt{-2}}$. 
 
 \smallskip	
 $b)$ We have $\Gamma $-equivariant isomorphisms
 \[H^{1,1}(C\times C)=H^2(C,\C)\oplus H^2(C,\C)\oplus (H^{1,0}\boxtimes H^{0,1})\oplus (H^{0,1}\boxtimes H^{1,0})=\C^2\oplus \End(H^0(C,K_C))^{{\scriptscriptstyle \oplus}2}\](where $\Gamma $ acts trivially on $\C^2$),
 hence $H^{1,1}((C\times C)/\Gamma )=\C^2\oplus \End_\Gamma (H^0(C,K_C))^{{\scriptscriptstyle \oplus}2}$. 
 
As a $\Gamma $-module we have $H^0(C,K_C)=L\oplus V$, where $\Gamma $ acts trivially on $L$ and $V$ is the sum of the 3 nontrivial one-dimensional representations
of $\Gamma $. Thus\[ \End_\Gamma (H^0(C,K_C)) = \mathbb{M}_2(\C)\times \C^3 \ .\]
Similarly we have $\NS((C\times C)/\Gamma )\otimes \Q=\Q^2 \oplus (\End_\Gamma (JC)\otimes \Q)$ and
\[ \End_\Gamma (JC)\otimes \Q = (\End (JH) \otimes \Q) \,\times  \, (\End_\Gamma (JF)\otimes \Q)^3 =\mathbb{M}_2(\Q(\sqrt{-2}))\times  \Q(i)^3\ ,\]
hence the result.\qed
 
 \begin{cor}[\cite{ST}]\label{ST}
Let $\Sigma \subset \P^6$ be the \emph{surface of cuboids}, defined by
\[ t^2= x^2+y^2+z^2 \quad , \quad u^2= y^2+z^2 \quad , \quad v^2=x^2+z^2 \qquad w^2=   x^2+y^2 \ .\] 
$\Sigma $ has $48$ ordinary nodes;  its minimal desingularization $S$ is $\rho $-maximal.\end{cor}
Indeed $\Sigma $ is a quotient of $(C\times C)/\Gamma $ \cite{Bt}.\qed

\smallskip	
(The result has been obtained first in \cite{ST} with a very different method.)

\medskip	
\section{Other examples}
\subsection{Elliptic modular surfaces}\label{mod}
Let $\Gamma $ be a finite index subgroup of $\mathrm{SL}_2(\Z)$ such that $-I\notin \Gamma $. The group $\mathrm{SL}_2(\Z)$ acts on the Poincar\'e upper half-plane $\mathbb{H}$; let $\Delta _\Gamma $ be the compactification of the Riemann surface
$\mathbb{H}/\Gamma $. 
The universal elliptic curve over $\mathbb{H}$ descends to 
$\mathbb{H}/\Gamma $, and extends to a smooth projective surface $B_\Gamma $ over $\Delta _\Gamma $, the \emph{elliptic modular surface} attached to $\Gamma $.
In \cite{S0} Shioda proves that  $B_\Gamma $  \emph{is $\rho $-maximal}. \footnote{  I am indebted to I. Dolgachev and B. Totaro for pointing out this reference.}

Now take $\Gamma =\Gamma (5)$, the kernel of the reduction map $\mathrm{SL}_2(\Z)\rightarrow \mathrm{SL}_2(\Z/5)$. In \cite{Li}
 Livne  constructed a $\Z/5$-covering $X\rightarrow B_{\Gamma (5)}$, branched along the sum of the 25 5-torsion sections of $B_{\Gamma (5)}$ . The surface $X$ satisfies $c_1^2=3c_2\,(=225)$, hence it is a ball quotient and therefore rigid. By analyzing 
the action of $\Z/5$ on $H^{1,1}(X)$ Livne shows that $H^{1,1}(X)$ is not defined over $\Q$, hence $X$ \emph{is not $\rho $-maximal}. This seems to be  the only known example 
of a surface which cannot be deformed to a  $\rho $-maximal surface. 
\medskip	
\subsection{Surfaces with $p_g=K^2=1$}

The minimal surfaces with $p_g=K^2=1$ have been studied by Cata\-nese \cite{C} and Todorov \cite{T1}. Their canonical model 
is a complete intersection of type $(6,6)$  in the weighted projective space $\P(1,2,2,3,3)$. The moduli space $\mathcal{M}$ is smooth of dimension 18.
\begin{prop}
The $\rho $-maximal surfaces are dense in $\mathcal{M}$.
\end{prop}
\pr 
We can replace $\mathcal{M}$ by the Zariski open subset $\mathcal{M}_{a}$ parametrizing  surfaces  with ample canonical bundle. Let $S\in\mathcal{M}_{a}$, and let
 $f:\mathcal{S}\rightarrow (B,\mathrm{o})$ be a local versal deformation of $S$, so that  $S\cong\mathcal{S}_{\mathrm{o}}$. 
Let $L$ be the lattice $H^2(S,\Z)$, and $k\in L$ the class of $K_S$. We may assume that $B$ is simply connected and fix an isomorphism of local systems $R^2f_*(\Z)\iso L_B$, compatible with the cup-product and mapping the canonical class $[K_{\mathcal{S}/B}]$ onto $k$. This induces for each $b\in B$ an isometry $\varphi _b: H^2(\mathcal{S}_b,\C)\iso L_{\C}$, which maps $H^{2,0}(\mathcal{S}_b)$ onto a line  in $L_{\C}$; the corresponding point $\wp(b)$ of $\P(L_{\C})$ is the period of $\mathcal{S}_b$. It belongs to the complex manifold
\[\Omega := \{[x]\in \P(L_{\C})\ |\ x^2=0\ ,\ x.k=0\ ,\ x.\bar{x}>0\}\ .\]
Associating to $x\in\Omega $ the real 2-plane $P_x:=\lr{\mathrm{Re}(x), \mathrm{Im}(x)}\subset L_{\R}$ defines an isomorphism of $\Omega $ onto  the Grassmannian of positive oriented 2-planes in $L_{\R}$. 

The key point is that
the image of the \emph{period map} $\wp: B\rightarrow \Omega $ is open\cite{C}.
Thus we can find $b$ arbitrarily close to $\mathrm{o}$ such that the 2-plane $P_b$ is defined over $\Q$,   
 hence  $H^{2,0}(\mathcal{S}_b)\oplus H^{0,2}(\mathcal{S}_b)=P_b\otimes _{\R}\C$ is defined over $\Q$.\qed

\medskip	
\begin{rem}
The proof applies to all surfaces with $p_g=1$ for which the image of the period map is open (for instance to K3 surfaces); unfortunately this seems to be a rather exceptional situation.
\end{rem}

\subsection{Todorov surfaces}
In \cite{T2} Todorov constructed a series of regular surfaces with $p_g=1$, \break $ 2\leq  K^2\leq 8$, which provide counter-examples to the Torelli theorem. The construction is as follows :
 let $ K \subset \P^3$ be a Kummer surface. We choose $ k $ double points of $ K $ in general position (this can be done with $ 0\leq k\leq 6$), and a general quadric  $Q\subset \P^3 $ passing through these $ k $ points. The \emph{Todorov surface} $ S $ is the double covering of $ K $ branched along $ K\cap Q $ and the remaining $ 16-k $ double points. It is a minimal surface of general type with $ p_g=1 $, $ K^2= 8-k $, $q=0$. If moreover we choose $ K $ $ \rho $-maximal (that is, $ K=E^2/\{\pm 1\} $, where $ E $ is an elliptic curve with complex multiplication), then $ S $ is $ \rho $-maximal by Proposition \ref{map}.$b)$.

Note that by varying the quadric $ Q $ we get a continuous, non-constant family of $ \rho $-maximal surfaces. 

\subsection{Double covers}\label{P}

In \cite{P} Persson constructs $\rho $-maximal double covers of certain rational surfaces by allowing  the branch curve to adquire some simple singularities (see also \cite{BE}). He applies this method to find $\rho $-maximal surfaces in the following families :

$\bullet$ Horikawa surfaces, that is, surfaces on the ``Noether line'' $K^2=2p_g-4$, for $p_g\not\equiv -1$ $(\mathrm{mod.}\ 6)$;

$\bullet$ Regular elliptic surfaces;

$\bullet$ Double coverings of $\P^2$.

In the latter case the double plane admits (many) rational singularities; it is unknown whether there exists a $\rho $-maximal surface $S$ which is a double covering of $\P^2$ branched along a smooth curve of even degree $\geq 8$.

\subsection{Hypersurfaces and complete intersections}

Probably the most natural  families to look at are smooth surfaces in $\P^3$, or more generally complete intersections.  
Here we may ask for a smooth surface $S$, or for the minimal resolution of a surface with rational double points (or even any surface  deformation equivalent to a complete intersection of given type). Here are the  examples that we know of :  

$\bullet$ The quintic surface $x^3yz+y^3zt+z^3tx+t^3xy=0$ has four $A_9$ singularities; its minimal resolution is $\rho $-maximal \cite{Sc}. It is not yet known whether there exists a
smooth $\rho $-maximal quintic surface.

$\bullet$ The Fermat sextic is $\rho $-maximal (\S \ref{prod}, Corollary \ref{S6}).

$\bullet$ The  complete intersection $y^2-xz = v^2-uw= x^3+z^3+u^3+w^3=0$ of type $(2,2,3)$ in $\P^5$ has 12 nodes; its minimal desingularization is $\rho $-maximal (end of \S \ref{prod}).

$\bullet$ The surface of cuboids is a complete intersection  of type $(2,2,2,2)$ in $\P^6$ with 48  nodes; its minimal desingularization is $\rho $-maximal (\S \ref{cub}, Corollary \ref{ST}).

\medskip	
\section{The complex torus associated to a $\rho $-maximal variety}
For a $\rho $-maximal variety $X$,
let $T_X$ be the $\Z$-module $H^2(X,\Z)/\mathrm{NS}(X)$. We have a decomposition
\[T_X\otimes \C= H^{2,0}\oplus H^{0,2}\]
defining a weight 1 Hodge structure on $T_X$, hence a complex torus $\mathcal{T}:=H^{0,2}/p_2(T_X)$, where $p_2:T_X\otimes \C\rightarrow  H^{2,0}$ is the second projection. Via the isomorphism $H^{0,2}=H^2(X,\O_X)$, $\mathcal{T}_X$ is identified with the cokernel of the natural map $H^2(X,\Z)\rightarrow H^2(X,\O_X)$.

The exponential exact sequence gives rise to an exact sequence
\[0\rightarrow \mathrm{NS}(X)\longrightarrow H^2(X,\Z)\longrightarrow  H^2(X,\O_X)\longrightarrow H^2(X,\O_X^*)\qfl{\partial} H^3(X,\Z)\ ,\] 
hence to a short exact sequence
\[0\rightarrow \mathcal{T}_X\longrightarrow H^2(X,\O_X^*)\qfl{\partial} H^3(X,\Z)\]
so that $\mathcal{T}_X$ appears as the ``continuous part'' of the group $H^2(X,\O_X^*)$.

\medskip	
\begin{ex}
Consider the elliptic modular surface $B_\Gamma $ of section \ref{mod}. The space $H^0(B_\Gamma ,K_{B_\Gamma })$ can be identified with the space of cusp forms of weight 3 for $\Gamma $; then the torus $\mathcal{T}_{B_\Gamma }$ is the complex torus associated to this space by Shimura (see \cite{S0}).
\end{ex}
\begin{ex}
Let $X=C\times C'$, with $JC$ isogenous to $E^g$ and $JC'$ to $E^{g'}$ (Proposition \ref{cc'}).  The torus $\mathcal{T}_X$ is the cokernel of the  map
 \[i\otimes i':H^1(C,\Z)\otimes H^1(C',\Z)\rightarrow H^1(C,\O_C)\otimes H^1(C',\O_{C'})\]
 where $i$ and $i'$ are the embeddings $H^1(C,\Z)\mono H^1(C,\O_C)$ and $H^1(C',\Z)\mono H^1(C',\O_{C'})$.
 
 We want to compute $\mathcal{T}_X$ up to isogeny, so we may replace the left hand side by a finite index sublattice. Thus, writing $E=\C/\Gamma $,
 we may identify $i$  with the diagonal embedding $\Gamma ^g\mono \C^g$, and similarly 
 for $i'$; therefore $i\otimes i'$ is the diagonal embedding of $(\Gamma \otimes \Gamma )^{gg'}$ in $\C^{gg'}$. Put $\Gamma =\Z+\Z\tau $;
 the image $\Gamma '$ of $\Gamma \otimes \Gamma $ in $\C$ is spanned by $1,\tau ,\tau ^2$; since $E$ has complex multiplication, $\tau $  is a quadratic number, hence $\Gamma $ has finite index in $\Gamma '$. Finally we obtain that $\mathcal{T}_X $\emph{is isogenous to} $E^{gg'}$. 
 
\end{ex}

For the surface $X=(C\times C)/\Gamma $ studied in \S \ref{cub} an analogous argument shows that $\mathcal{T}_X$ is isogenous to $A=E_i^4\times E_{\sqrt{-2}}^3$. This is still an abelian variety of type CM, in the sense that $\End(A)\otimes \Q$ contains an \'etale $\Q$-algebra of maximal dimension $2\dim(A)$. There  seems to be no reason why this should hold in general. However it is true in the special case 
 $h^{2,0}=1$ (e.g. for holomorphic symplectic manifolds):

\begin{prop}
If $h^{2,0}(X)=1$, the torus $\mathcal{T}_X$ is an elliptic curve with complex multiplication.
\end{prop}
\pr  Let $T'_X$ be the pull back of $H^{2,0}+H^{0,2}$ in $H^2(X,\Z)$;
then $p_2(T'_X)$ is a sublattice of finite index in $p_2(T_X)$. 
Choosing an ample class $h\in H^2(X,\Z)$ defines a
 quadratic form on $H^2(X,\Z)$ which is positive definite on $T'_X$.
Replacing again $T'_X$ by a finite index sublattice we may assume that it admits an orthogonal basis $(e,f)$ with $e^2=a$, $f^2=b$. Then $H^{2,0}$ and $H^{0,2}$ are the two 
isotropic lines of $T'_X\otimes \C$; they are spanned by the vectors $\omega =e+\tau f$ and $\bar{\omega }=e-\tau f$, with $\tau ^2=-a/b$. We have $e=\frac{1}{2} (\omega +\bar{\omega }) $ and $f=\frac{1}{2\tau  } (\omega -\bar{\omega }) $; therefore  multiplication by $\frac{1}{2\tau  } \bar{\omega}  $ induces an isomorphism of $\C/(\Z+\Z \tau  )$ onto $H^{0,2}/p_2(T'_X)$, hence $\mathcal{T}_X$ is isogenous to $\C/(\Z+\Z \tau  )$ and $\End(\mathcal{T}_X)\otimes \Q=\Q(\tau  )=\Q(\sqrt{-\mathrm{disc}(T'_X) })$.\qed

\medskip
\section{Higher codimension cycles}

A natural generalization of the question considered here is to look for varieties $X$ for which the group  
$ H^{2p}(X,\Z)_{\mathrm{alg}}$  of algebraic classes in $H^{2p}(X,\Z)$ has maximal rank $h^{p,p}$. Very few nontrivial cases seem to be known. The following is essentially due to Shioda :
\begin{prop}
Let $F^n_d$ be the Fermat hypersurface of  degree $d$ and even dimension $n=2\nu $.  For $d=3,4$, the group $ H^{n}(F^n_d,\Z)_{\mathrm{alg}}$   has maximal rank $h^{\nu ,\nu }$. 
\end{prop}
\pr According to \cite{S1} we have
\[\rk  H^{n}(F^n_3,\Z)_{\mathrm{alg}}=1+\frac{n!}{(\nu  )!^2} \qquad  \mbox{and}  
\qquad \rk  H^{n}(F^n_4,\Z)_{\mathrm{alg}}= \sum_{k=0}^{k=\nu +1 } \dfrac{(n+2)!}{(k!)^2(n+2-2k)!}\ \cdot \] 
 On the other hand,  let $R^n_d:= \C[X_0,\ldots ,X_{n+1}]/(X_0^{d-1},\ldots , X_{n+1}^{d-1})$ be the Jacobian ring of $F^n_d$; Griffiths theory \cite{G} provides an isomorphism of the primitive cohomology $H^{\nu ,\nu }(F^n_d)_{\mathrm{o}}$ with the component of degree $(\nu +1)(d-2)$ of $R^n_d$. Since this ring is the tensor product of $(n+2)$ copies of $\C[T]/(T^{d-1})$, its Poincar\'e series $\sum_k \dim (R^n_d)_kT^k$ is $(1+T+\ldots +T^{d-2})^{n+2}$. Then an elementary computation gives the result.\qed 
 
 \medskip	
In the particular case of cubic fourfolds we have more examples :
\begin{prop}
Let $F$ be a cubic form in 3 variables, such that the curve $F(x,y,z)=0$ in $\P^2$ is an elliptic curve with complex multiplication; let $X$ be the cubic fourfold  defined by 
$F(x,y,z)+F(u,v,w)=0$ in $\P^5$. The group  $ H^4(X,\Z)_{\mathrm{alg}}$ has maximal rank $h^{2,2}(X)$. 
\end{prop}
\pr 
Let $u$ be  the  automorphism of $X$ defined by $u(x,y,z;u,v,w)=(x,y,z;\omega u,\omega v,\omega w)$. We observe that $u$ acts trivially on the (one-dimensional) space $H^{3,1}(X)$. Indeed Griffiths theory \cite{G} provides a canonical isomorphism 
\[\mathrm{Res}:\ H^0(\P^5,K_{\P^5}(2X))\iso H^{3,1}(X)\ ;\]
the space $H^0(\P^5,K_{\P^5}(2X))$ is generated by the meromorphic form $\dfrac{\Omega }{G^2}$ ,  with \[ \Omega =xdy\wedge dz\wedge du \wedge dv \wedge dw - y dx \wedge dz\wedge du \wedge dv \wedge dw + \ldots \quad \mbox{and} \quad G=F(x,y,z)+F(u,v,w)\ .\] The automorphism $u$ acts  trivially on this form, and therefore on $H^{3,1}(X)$.

Let $F$ be the variety of lines contained in $X$. We recall from \cite{BD} that $F$ is a holomorphic symplectic fourfold, and that there is a natural isomorphism of Hodge structures $\alpha  :H^4(X,\Z)\iso H^2(F,\Z)$. Therefore the automorphism $u_F$ of $F$ induced by $u$ is symplectic.
Let us describe its fixed locus. 

The  fixed locus of $u$ in $X$ is the union of the plane cubics $E$ given by $ x=y=z=0$ and $E'$ given by $ u=v=w=0$. A line in $X$ preserved by $u$ must have (at least) 2 fixed points, hence must meet both $E$ and $E'$; conversely, any line joining a point of $E$ to a point of $E'$ is contained in $X$, and preserved by $u$. This identifies  the fixed locus $A$ of $u_F$ to the abelian surface $E\times E'$. Since $u_F$ is symplectic $A$ is a symplectic submanifold, that is, the restriction map $H^{2,0}(F)\rightarrow H^{2,0}(A)$ is an isomorphism. By our hypothesis $A$ is $\rho $-maximal, so $F$ is $\rho $-maximal by Proposition \ref{map}. Since $\alpha $ maps $H^4(X,\Z)_{\mathrm{alg}}$ onto $NS(F)$ 
this implies the Proposition. \qed

\end{document}